\documentclass[3p,preprint]{elsarticle}

\usepackage{amsmath}
\usepackage{amsfonts}
\usepackage{amssymb,latexsym}
\usepackage{amsthm}
\usepackage{graphicx}

\newtheorem{Remark}{Remark}[section]

\usepackage{color}
\usepackage{newlfont}
\usepackage{subfigure}
\usepackage{verbatim}
\usepackage{rotating}
\usepackage{multirow}
\usepackage{units}
\usepackage{bm}
\usepackage[english]{babel}
\usepackage[utf8]{inputenc}
\usepackage{algorithm}
\usepackage[noend]{algpseudocode}
\usepackage{booktabs}
\usepackage{caption}
\usepackage{hyperref}

\newcommand{\red}[1]{{\color{red} #1}}

\journal{a journal}

\begin{document}

\begin{frontmatter}

\title{inRS: implementing the indicator function \\ 
of NURBS-shaped planar domains}



\author[address-PD,address-GNCS]{Alvise Sommariva}
\ead{alvise@math.unipd.it}

\author[address-PD,address-GNCS]{Marco Vianello}
\ead{marcov@math.unipd.it}

\cortext[corrauthor]{Corresponding author}

\address[address-PD]{University of Padova, Italy}

\address[address-GNCS]{Member of the INdAM Research group GNCS}

\begin{abstract}

We provide an algorithm that implements the indicator function of   
NURBS-shaped planar domains, tailored to the fast computation on 
huge point clouds, together with the corresponding Matlab code. 
\end{abstract}



\end{frontmatter}


\section{Introduction}

NURBS-shaped domains produced by CAGD algorithms play a central
role in digital design and modelling processes. The capability
of locating
quasi-uniform or random sample points in such domains 
can be useful in a vast range of applications, for example within 
several meshfree bivariate approximation
algorithms developed in the last twenty years, among which we may 
quote (without any pretence of completeness) 
kernel-based 
and partition-of-unity collocation methods \cite{CB15,FMC15}, 
construction of algebraic 
cubature formulas \cite{GUND21,SV21-first,SV21} potentially useful for curved 
FEM/VEM elements \cite{BdVRV19,SFM11}, 
compressed MC/QMC integration \cite{BDME16,H21}, compressed polynomial 
regression \cite{BPV20,PSV17}.

Though the efficient computation of the indicator function of 
general linear polygons has received much attention in the algorithmic 
literature and deserved sophisticated implementations (cf. e.g. 
\cite{HA01}), such as 
the Matlab {\sc inpolygon} function (and the {\sc isinterior} 
function in the {\sc polyshape} envinronment), or the recent 
{\sc inpoly2} function \cite{E21}, the same cannot be 
apparently said concerning NURBS-shaped curved polygons (at least in Matlab), 
despite their 
relevance in applications.     

Extending and improving an approach already explored in 
\cite{SV21-first,SV21}, as well as resorting to some ideas used 
in \cite{E21} for linear polygons,  
we try to fill this lack by providing an efficient and robust 
Matlab implementation 
of the indicator function of NURBS-shaped Jordan domains, 
based on the topological notions of crossing number 
(even-odd rule) and winding number. 
A key tool is encapsulating the domain boundary by a finite number of 
Cartesian rectangles, in all of which it is the graph 
of a local monotone Cartesian function. We have tried to optimize 
all the algorithm 
blocks and to conveniently manage the critical situations,  
the present strategy being mainly tailored to fast computation 
of the indicator function on huge point clouds. 
The corresponding  
Matlab code \cite{SV22-soft}, that could be useful in many design
and modelling applications, is freely available to the 
scientific community.   

\section{NURBS-shaped indicator function}

In this section we discuss an algorithm for the computation of the 
indicator function of bidimensional NURBS-shaped Jordan domains, based 
on the topological notions of 
crossing number and winding number, 
whose main lines have already appeared in \cite{SV21}. The present 
implementation is more efficient and robust, and for the reader's 
convenience we explain in some detail the whole construction.     

Consider a Jordan domain $\Omega \subset {\mathbb{R}}^2$,
such that 
\begin{itemize}
\item its boundary $\partial \Omega$ is a simple curve 
described by parametric equations 
\begin{equation} \label{curve}
\Gamma(t)=(\alpha(t),\beta(t))\;,\;\;t \in I=[t_{min},t_{max}]\;,
\end{equation}
where ${\alpha},{\beta} \in C(I)$,
$\Gamma(t_{min})=\Gamma(t_{max})$;
\item there are partitions $\{I_k\}_{k=1,\ldots,K}$ of $I$,
and $\{I_{k,j}\}_{j=1,\ldots,m_k}$ of each $I_k= 
[t_k, t_{k+1}]$,
such that the restrictions of ${\alpha,\beta}$
to each $I_k$ are {\sl{rational splines}},
w.r.t. the subintervals $\{I_{k,j}\}$.
\end{itemize}

We shall denote the local rational splines as 
\begin{equation} \label{Ik}
{\alpha}(t)=\frac{u_k(t)}{v_k(t)}, \,\,\,
{\beta}(t)=\frac{w_k(t)}{z_k(t)}, \,\,\, t \in I_k\;,
\end{equation}
where the numerators $u_k$, $w_k$ as well as the denominators 
$v_k$, $z_k$ are polynomial splines on $I_k$,
sharing the same knots and having degree, respectively,
$p_k$ and $q_k$.
Notice that since ${\alpha,\beta}$ are globally continuous,
the denominators $v_k$, $z_k$ do not vanish 
in the closed interval $I_k$. Moreover, we have that 
$\partial \Omega=\cup_{k=1}^{M} 
(V_k \frown V_{k+1})$, with the convention that $V_{K+1}=V_1$, 
where the {\em vertices} $\{V_k\}$ 
can be corner points or even cusps of the boundary. 
On the other hand, in each subinterval $I_{k,j}\subset I_k$ we have that 
\begin{equation} \label{Ikj}
{\alpha}(t)=\frac{u_{k,j}(t)}{v_{k,j}(t)}, \,\,\,
{\beta}(t)=\frac{w_{k,j}(t)}{z_{k,j}(t)}, \,\,\, t \in I_{k,j}\;,
\end{equation}
where the numerators $u_{k,j},w_{k,j}$ are polynomials of degree $p_k$ 
and the denominators $v_{k,j},z_{k,j}$ are polynomials of degree $q_k$.

We are particularly interested in the case when
the boundary is locally a $p$-th degree NURBS curve
{\cite[p.117]{PIEGL}}, i.e. the curvilinear side
$V_k \frown V_{k+1}$ has the following parametrization  

\begin{equation} \label{side}
\Gamma(t)=\frac{\sum_{i=1}^{m_k} B_{i,p}(t) \,\lambda_{i,k} C_{i,k} }
{\sum_{i=1}^{m_k} B_{i,p}(t) \,\lambda_{i,k}}, \,\,\, t 
\in [t_k, t_{k+1}]\;,
\end{equation}

\noindent where $\{C_{i,k}\}_{i=1}^{m_k}\subset \mathbb{R}^2$ are the
{\it{control points}},  $\{\lambda_{i,k}\}_{i=1}^{m_k}$ are the weights
and $\{B_{i,p}\}_{i=1}^{m_k}$ are the $p$-th degree B-spline basis
functions {\cite[p.87]{DEBOOR}} defined on 
a suitable knot vector.
We stress that the case of standard polynomial splines 
(already treated in \cite{SV21-first}) is included 
in (\ref{curve})-(\ref{Ikj}), 
and that the whole construction below can be extended  
also to other domains with Rational Spline (RS) boundary, 
such as rational Bezier curves. 

\begin{figure}[!htbp]
  \centering
   {\includegraphics[scale=0.2,clip]{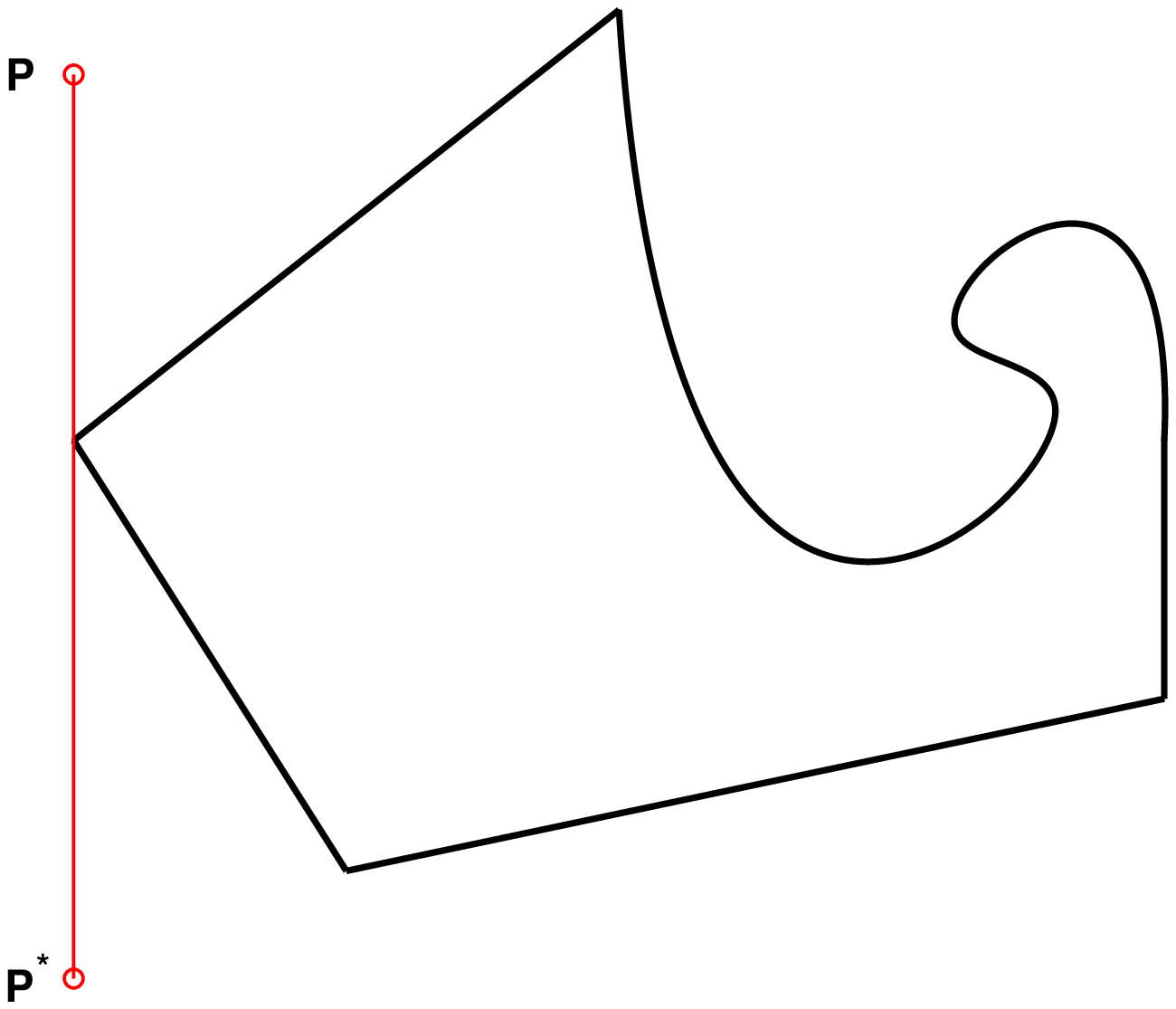}}
   {\includegraphics[scale=0.2,clip]{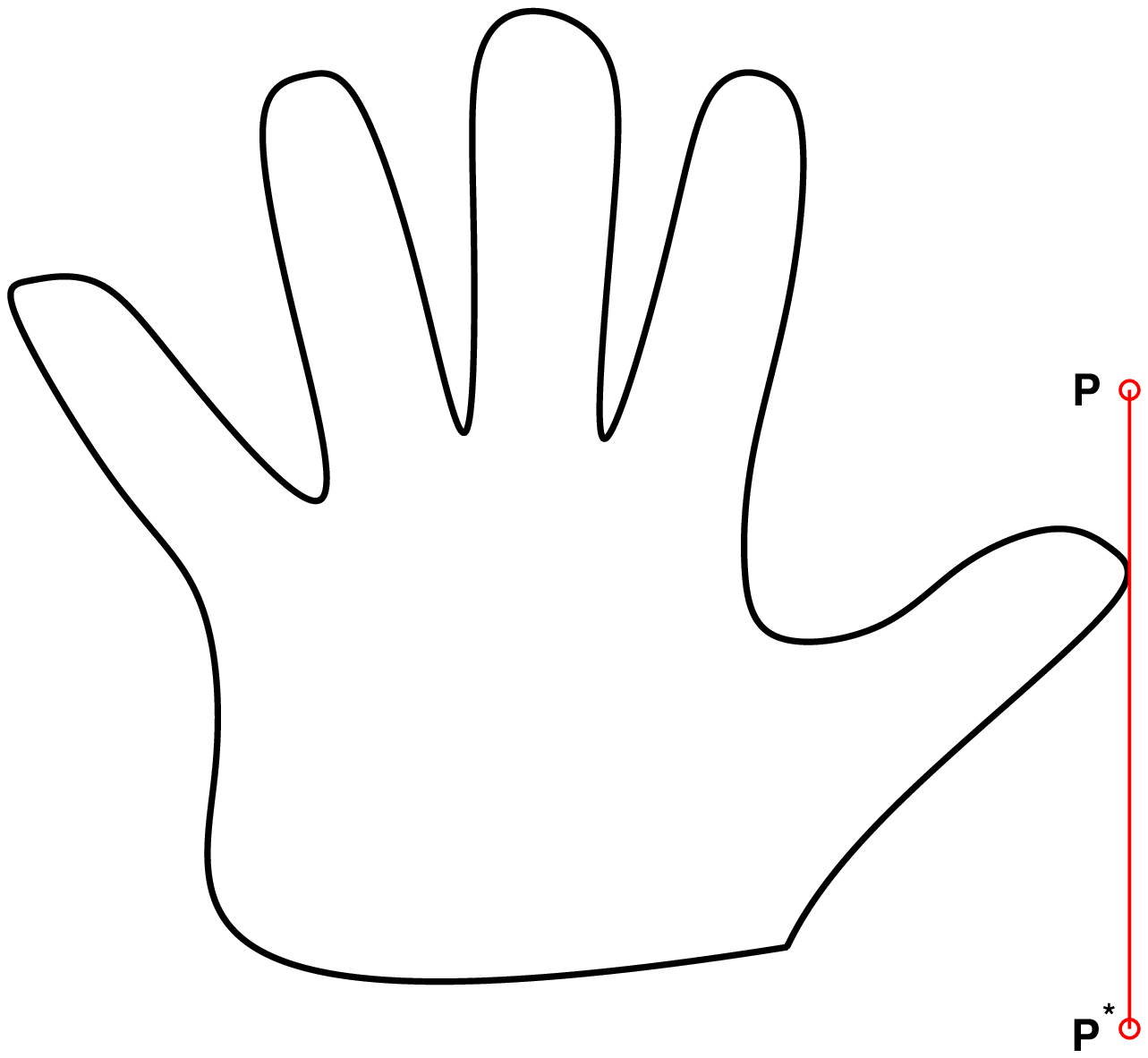}}
   {\includegraphics[scale=0.2,clip]{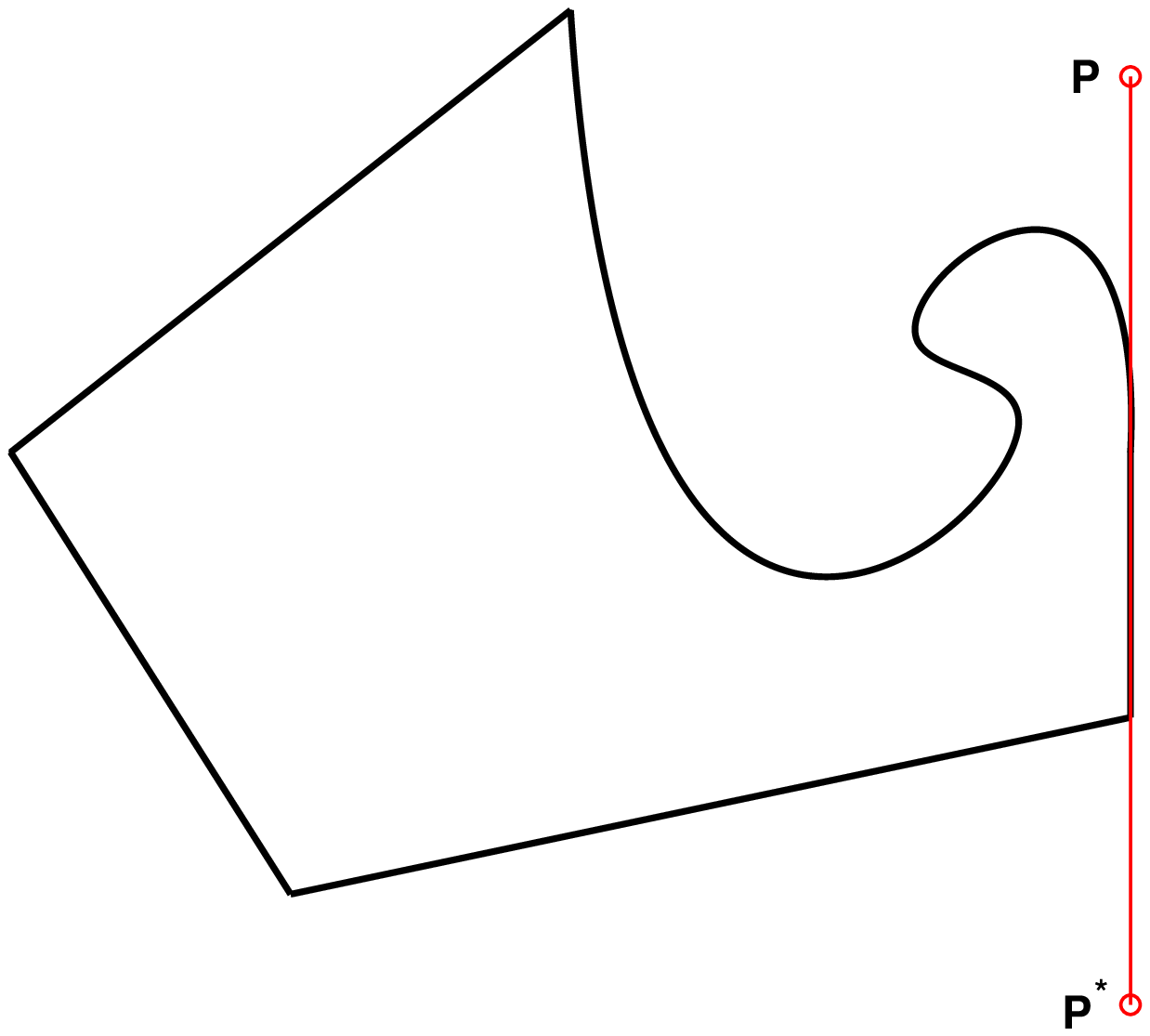}}
 \caption{Critical situations for the crossing number
on curvilinear domains.}
 \label{fig_1CRI}
 \end{figure}

We base our indicator function algorithm $inRS$  
on the well-known {\it{Jordan curve
theorem}}, which implies that a point $P$ belongs to a Jordan domain 
$\Omega$ if and only if,
taking a point $P^* \notin \Omega$, the segment
${\overline{P^*P}}$ {\it crosses} $\partial \Omega$
an odd number of times; cf. e.g. \cite{J1893} and 
the nice paper \cite{H07} with the references therein.
There may be some ``critical'' cases, for instance when
${\overline{P^*P}}$ touches a vertex without crossing the boundary, 
or when it is tangent  
to the boundary; cf. Fig. \ref{fig_1CRI} where vertical segments 
are used as in our main implementation. In these cases the crossing
number strategy cannot be directly applied and alternatives have 
to be adopted, such as computation of the winding number.    


The first step of the algorithm consists in covering 
$\partial \Omega$ with a finite union
of suitable Cartesian rectangles, in all of which 
the boundary is the graph
of a local monotone Cartesian function. Each rectangle contains a portion 
of $\partial \Omega$ that
is parametrized by two {\it{rational functions}},
i.e. locally $({\alpha}(t),{\beta}(t))$ are the
ratio of two polynomials.
Once these covering rectangles have been determined, 
evaluating the crossing number $cross(P)$ (i.e. computing the indicator 
function at a given point $P$) becomes easy, requiring at most 
the solution of some polynomial equations.

To this purpose, first we compute in each $I_{k,j}$ the possible 
zeros of $\alpha'(t)=(u'_{k,j}v_{k,j}-u_{k,j}v'_{k,j})/v^2_{k,j}$ 
and $\beta'(t)=(w'_{k,j}z_{k,j}-w_{k,j}z'_{k,j})/z^2_{k,j}$, 
obtaining a finer partition, say $I_{k,j,s}$, 
with such zeros as endpoints. Notice that this requires solving 
in each $I_{k,j}$ the two polynomial equations of degree $p_k+q_k-1$ 
\begin{equation} \label{zerosalpha'beta'}
u'_{k,j}(t)v_{k,j}(t)-u_{k,j}(t)v'_{k,j}(t)=0\;,\;\;
w'_{k,j}(t)z_{k,j}(t)-w_{k,j}(t)z'_{k,j}(t)=0\;,
\end{equation}
which can be conveniently done in Matlab by the {\sc roots} function, 
that comptes the eigenvalues of the companion matrix. 

\begin{figure}[!htbp]
  \centering
   {\includegraphics[scale=0.2,clip]{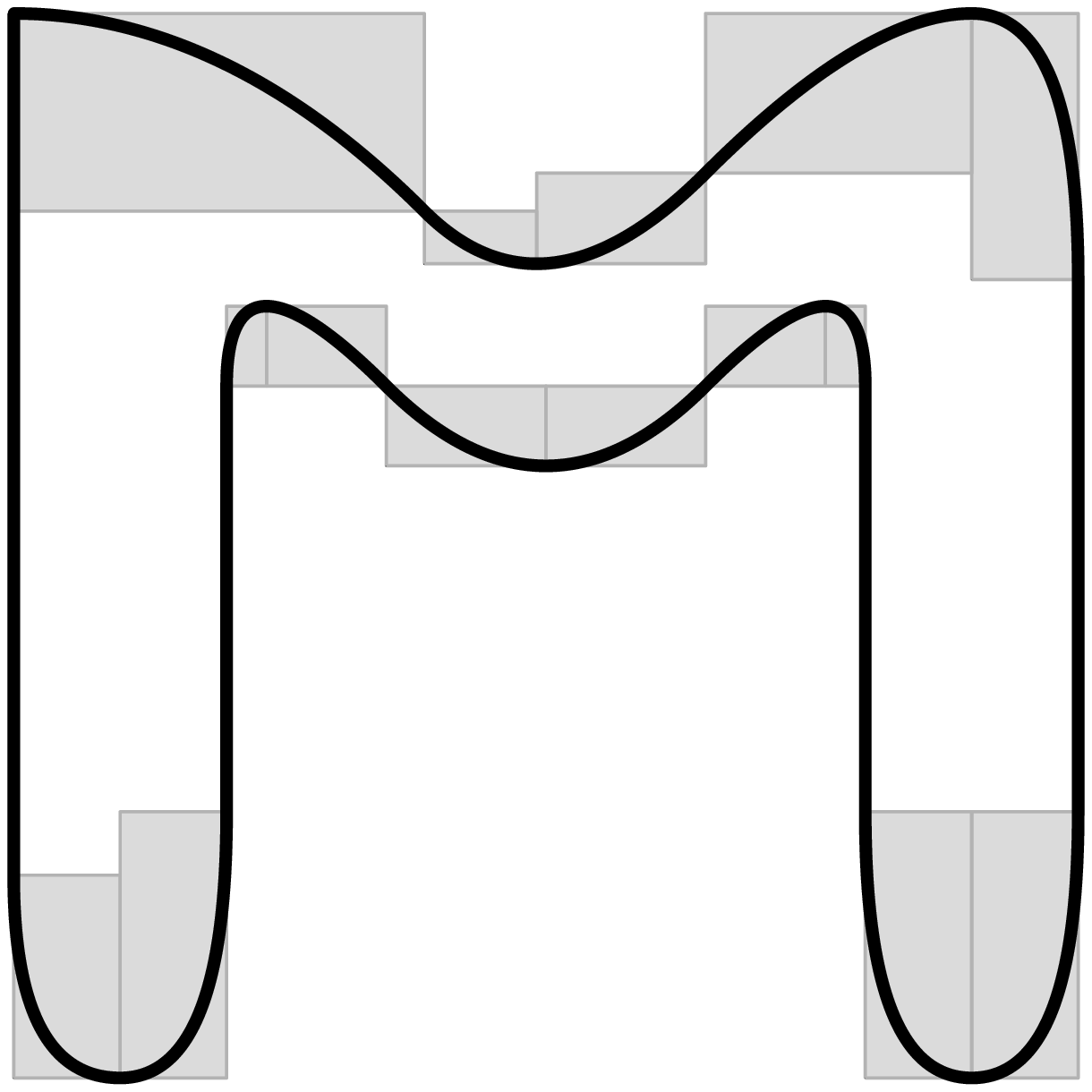}}
     {\includegraphics[scale=0.2,clip]{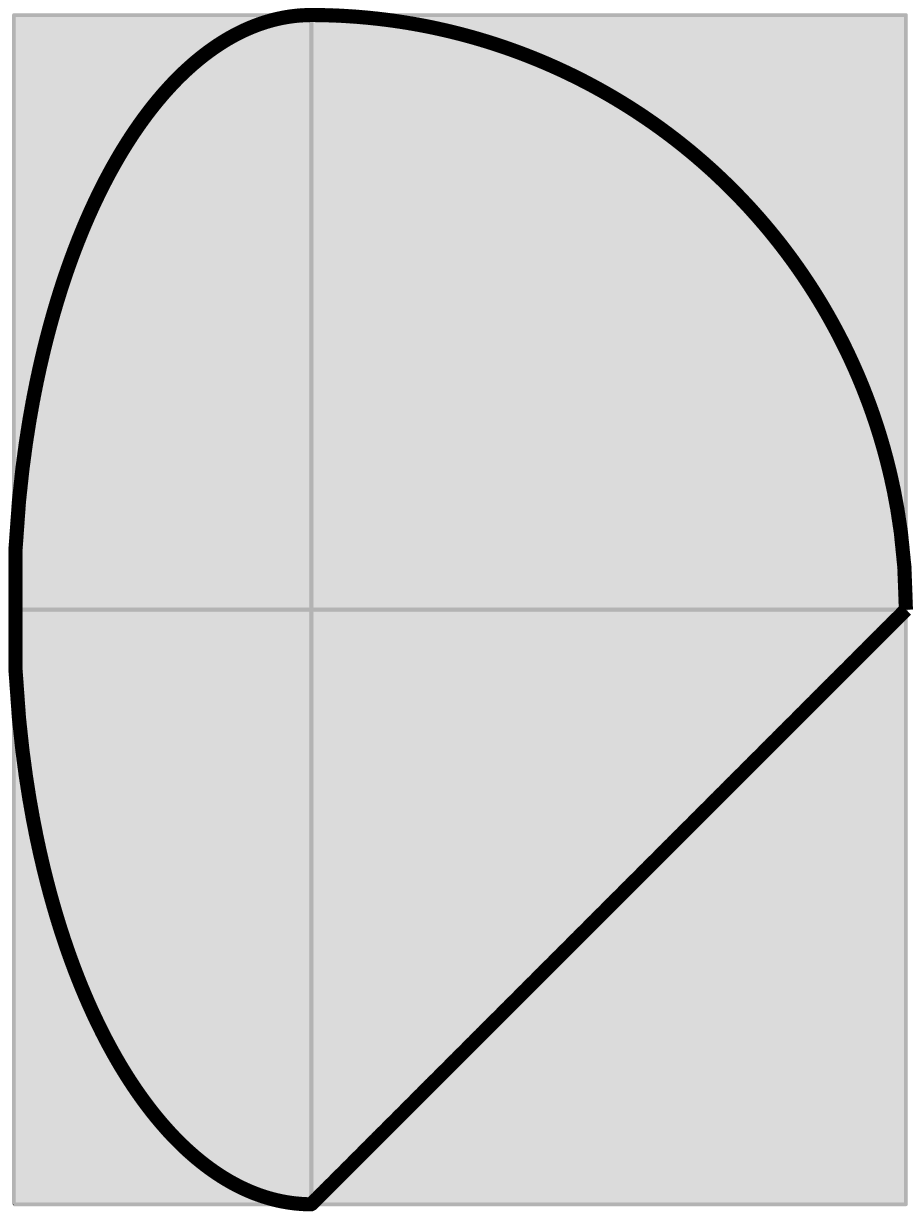}}
       {\includegraphics[scale=0.2,clip]{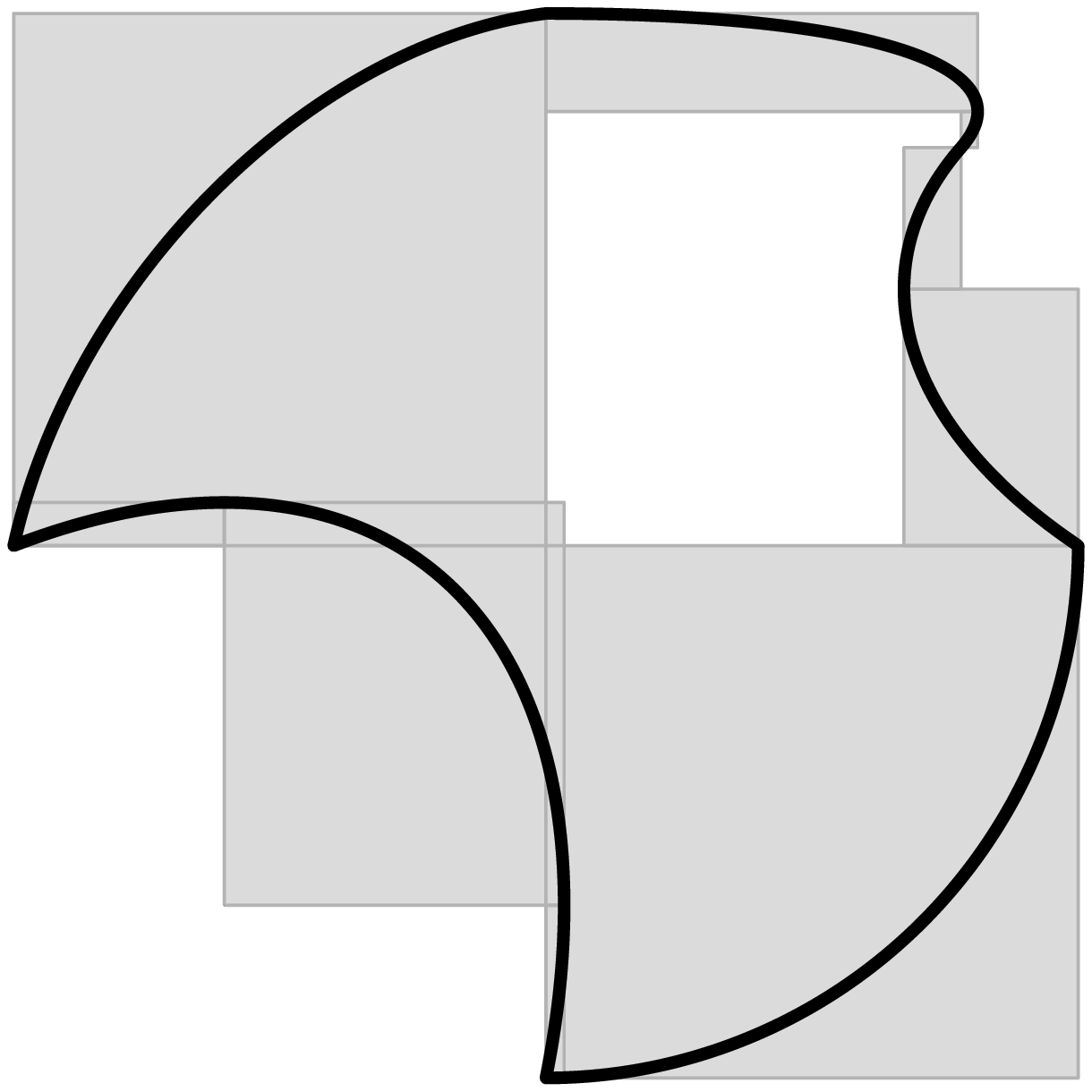}}
 \caption{Three NURBS-shaped domains
whose boundary may contain arc of circles, ellipses, segments as well
as other NURBS blocks, together with the corresponding monotone boxes.
Notice that the latter may degenerate into a segment
(left figure) and may overlap
(right figure).}
 \label{fig2}
 \end{figure}

Now, $\alpha(t)$ and $\beta(t)$ being polynomials, 
if not constant are strictly monotone in each $I_{k,j,s}$, 
and thus the boundary curve is there 
the graph of a strictly monotone Cartesian function, 
with local bounding box 
\begin{equation} \label{mbox}
\mathcal{B}_{k,j,s}=\left[a_{k,j,s},b_{k,j,s}\right]\times 
\left[c_{k,j,s},d_{k,j,s}\right]
\end{equation}
$$
a_{k,j,s}=\min_{t\in I_{k,j,s}}\alpha(t)\;,\;
b_{k,j,s}=\max_{t\in I_{k,j,s}}\alpha(t)
$$
$$
c_{k,j,s}=\min_{t\in I_{k,j,s}}\beta(t)\;,\;
d_{k,j,s}=\max_{t\in I_{k,j,s}}\beta(t)
$$
that we call a {\em monotone box}. The local minima and maxima
can be determined via the values of $\alpha,\beta$ at the endpoints 
of $I_{k,j,s}$. The case of $\alpha(t)$ or $\beta(t)$ constant are treated 
separately and correspond to bounding boxes degenerating into a vertical 
or horizontal segment, respectively. 
Clearly, the union of such monotone boxes, that may overlap, covers
the whole boundary $\partial \Omega$; cf. Fig. \ref{fig2}. 
Moreover, a global bounding box for the whole $\Omega$ can 
be immediately determined by the upmost, downmost, leftmost and rightmost 
monotone boxes. 

Now, take a point $P=(\overline{x},\overline{y})\in \mathbb{R}^2$. 
If $P$ is out of the global bounding box, 
then clearly the indicator function at $P$ is null, $\chi_\Omega(P)=0$.
If $P$ is in the global bounding box, consider 
the vertical ``downward''
ray $\{x=\overline{x},y\leq \overline{y}\}$ and the monotone 
boxes $\mathcal{B}_\ell$ corresponding to the triples 
\begin{equation}
\label{whichboxes}
\{\ell=(k,j,s):\;a_{k,j,s}\leq \overline{x}\leq b_{k,j,s}\;,\;\;
\overline{y}\geq c_{k,j,s}\}\;, 
\end{equation}
that are the boxes that do intersect the ray and  
can be effectively determined by a fast vectorized search. Notice that 
if $\overline{y}>d_\ell$ (that is $P$ is over $\mathcal{B}_\ell$ ) 
the ray surely intersects 
at one point the boundary portion pertaining to such a monotone box, 
whereas if $P$ is in $\mathcal{B}_\ell$ 
the possible intersection can be ascertained by solving 
the polynomial equation of degree $\max{\{p_k,q_k\}}$ 
\begin{equation} \label{intersections}
u_{k,j}(t)-\overline{x}\,v_{k,j}(t)=0\;,\;\;t\in I_\ell\;,
\end{equation}
which again can be conveniently done in Matlab 
by the {\sc roots} function.
Indeed, let $t_\ell$ be the unique solution. Then if 
$\beta(t_\ell)=w_{k,j}(t_\ell)/z_{k,j}(t_\ell)\leq \overline{y}$ 
the ray intersects the boundary at $(\alpha(t_\ell),\beta(t_\ell))$ 
and the intersection 
counting must be increased by 1, whereas if 
$\beta(t_\ell)>\overline{y}$ it does not. 

At this point, let $cross(P)$ be the overall number of  
such intersections. If none of them is a critical point, 
$cross(P)$ is just the crossing number of Jordan curve theorem. 
This can be 
easily ascertained by checking that $\overline{x}$ is not the abscissa 
of a boundary point where $\alpha'$ vanishes ({\em point of vertical 
tangency}),  
or of a vertex $V_i=(\alpha(t_i),\beta(t_i))$ such that 
$\alpha'(t_i^-)\alpha'(t_i^+)<0$, i.e. where the boundary curve 
turns form left to right or conversely ({\em $x$-turning vertex}), 
both being possible intersections without crossing. 
Then $P$ belongs to $\Omega$ if and only if $cross(P)$ is odd, i.e. the 
indicator function at $P$ is 
\begin{equation} \label{chi1}
\chi_\Omega(P)=cross(P)\;\mbox{mod}\,2\;. 
\end{equation}
In addition, one may also know whether the point is on the boundary or 
in the interior of $\Omega$, by checking whether one of such intersections 
coincides with $P$ (up to a suitable numerical tolerance). 

Whenever critical points are present among the intersections, one may 
use the same procedure 
working with an horizontal ray. In the rare case 
of another 
failure still due to critical points, 
one can resort to a different topological index, i.e. the winding number 
$wind(P)\in \mathbb{Z}$ (cf. e.g. \cite{Kr99})  
\begin{equation} \label{winding}
wind(P)= \frac{1}{\mbox{length}(I)}\,\int_{I}
\frac {\beta'(t)\,(\alpha(t)-\overline{x})
-\alpha'(t)\,(\beta(t)-\overline{y})} { (\alpha(t)-\overline{x})^2
+  (\beta(t)-\overline{y})^2  }\,dt
\end{equation}
$$
=\frac{1}{\mbox{length}(I)}\,\sum_{k,j}\int_{I_{k,j}}
\frac {\beta'(t)\,(\alpha(t)-\overline{x})
-\alpha'(t)\,(\beta(t)-\overline{y})} { (\alpha(t)-\overline{x})^2
+  (\beta(t)-\overline{y})^2  }\,dt\;,\;\;
P=(\overline{x},\overline{y})
\notin \partial \Omega\;,
$$
since (assuming that the boundary is counterclockwise oriented) 
the indicator function at $P$ is   
\begin{equation} \label{chi2}
\chi_\Omega(P)=wind(P)\;. 
\end{equation}
We observe that the evaluation of the integral above, for 
example by high-precision Gaussian quadrature on each subinterval 
in view of the fact that $\alpha,\beta\in C^\infty(I_{k,j})$, 
can be difficult 
when $P$ is close to the boundary. In practice, however, there 
is no need to compute
such a quantity with a small error, recalling that $wind(P)$ 
is an integer (so that an error strictly less than $1/2$ would 
suffice). 

\begin{Remark}
{\em 
One may cover the boundary with a larger number of monotone boxes,  
by further partitioning each $I_{k,j,s}$ into a number of 
subintervals $I_{k,j,s,\tau}=[a_{k,j,s,\tau},b_{k,j,s,\tau}]$. 
In such a way the boxes become 
clearly smaller and we have a finer approximation of the NURBS-shaped 
boundary: see Fig. \ref{fig2}. Given a point 
$P=(\overline{x},\overline{y})$ the pertaining boxes 
correspond to the quadruples  
\begin{equation}
\label{whichboxes2}
\{\ell=(k,j,s,\tau):\;a_{k,j,s,\tau}\leq \overline{x}\leq b_{k,j,s,\tau}\;,\;\;
\overline{y}\geq c_{k,j,s,\tau}\}\;, 
\end{equation}
where $c_{k,j,s,\tau}=\min_{t\in I_{k,j,s,\tau}}\beta(t)=
\min \{\beta(a_{k,j,s,\tau}),\beta(b_{k,j,s,\tau})\}$.

On the other hand, the smaller the boxes the 
smaller the probability that a random point is inside some box,
and thus we can substantially reduce the number of equations 
like (\ref{intersections}) to be solved by {\sc roots} in order 
to find the intersections and speed-up 
the whole procedure, when a huge set of points has to be located. 
In practice, however, there is a suitable threshold for the number of 
sub-boxes, over which refining is no more convenient (such a threshold can 
roughly determined experimentally).    
\/}
\end{Remark}

\begin{figure}[!htbp]
  \centering
   {\includegraphics[scale=0.2,clip]{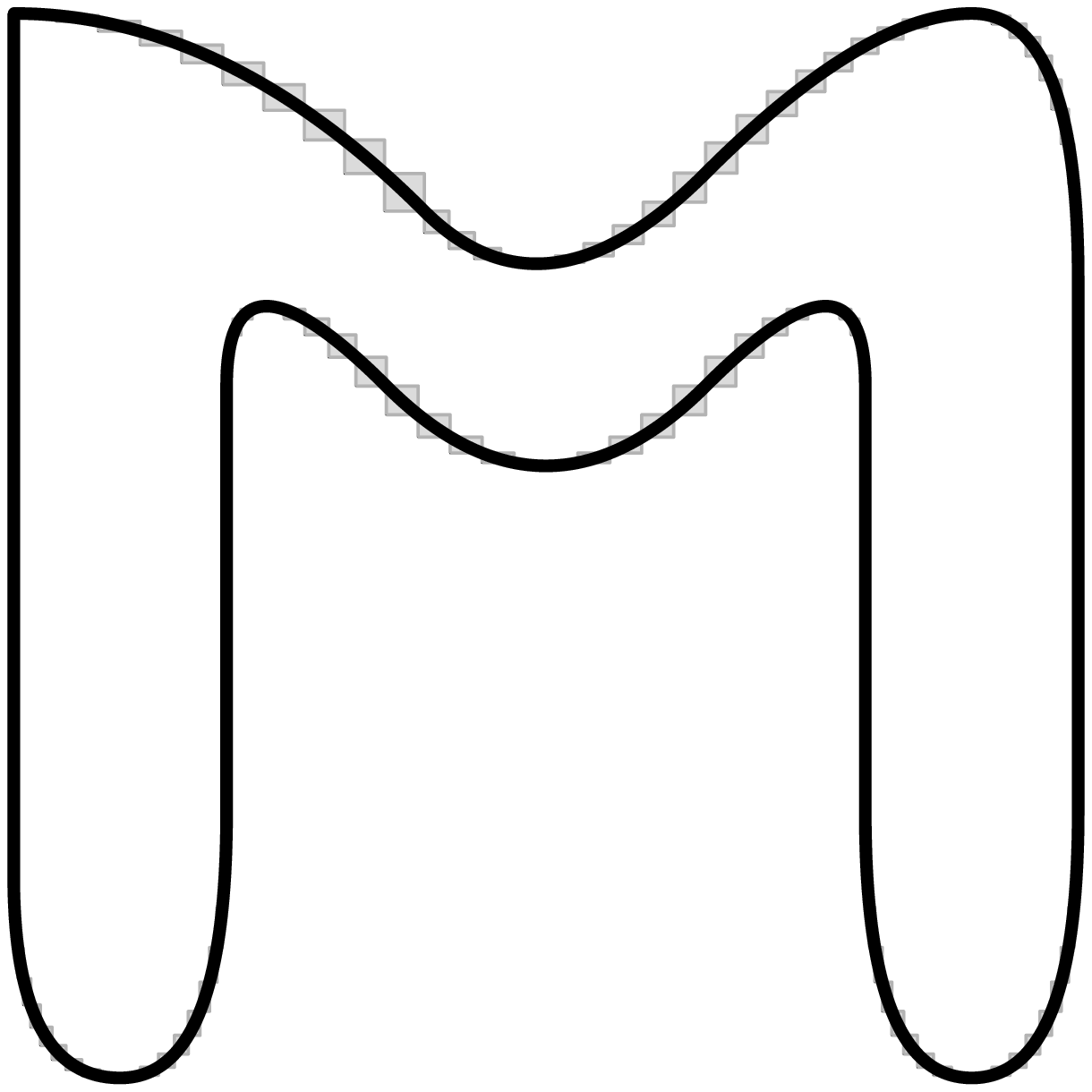}}
     {\includegraphics[scale=0.2,clip]{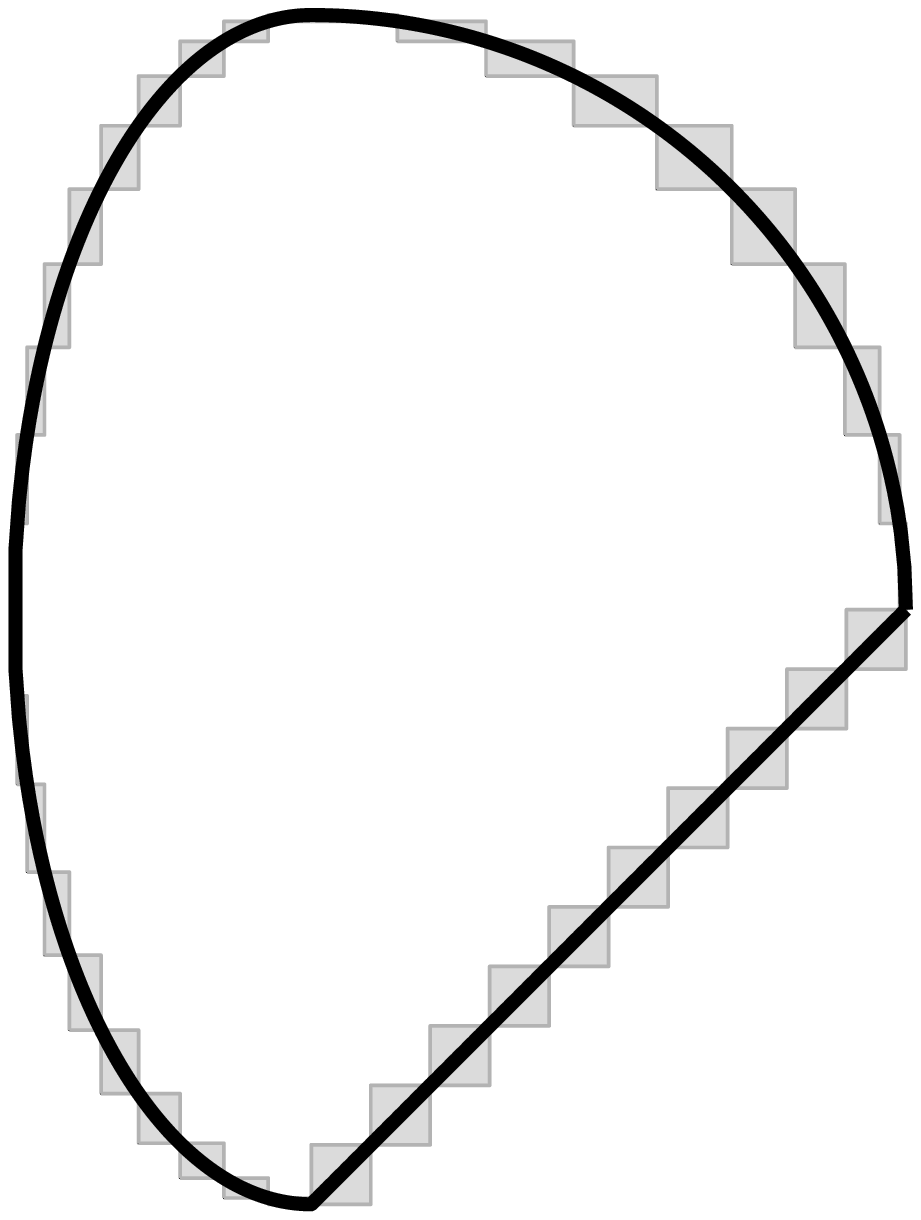}}
       {\includegraphics[scale=0.2,clip]{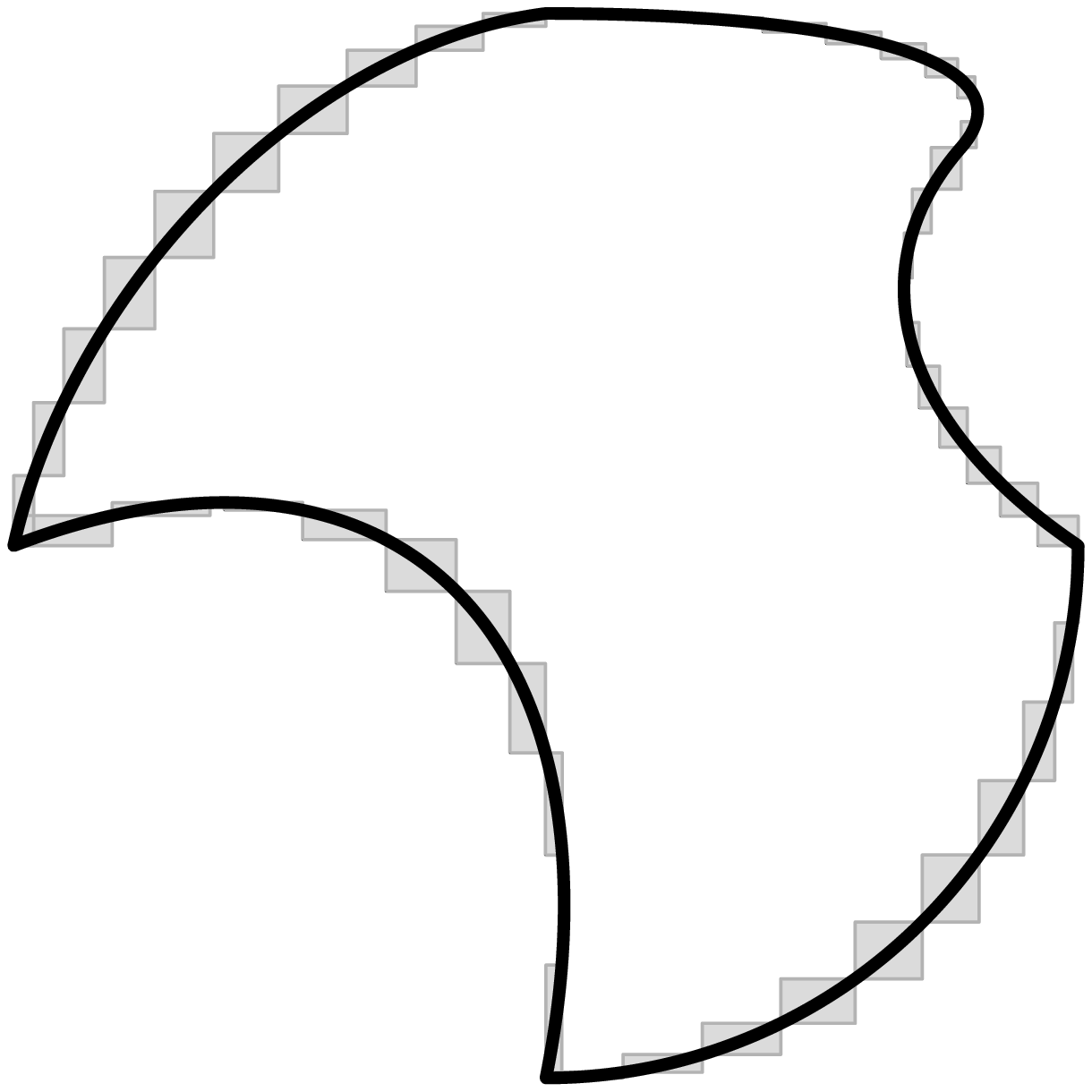}}
 \caption{Box refinement for the three NURBS-shaped domains 
of Fig. \ref{fig2}.}
 \label{fig3}
 \end{figure}

To fix ideas, the guidelines of the implementation 
based on the discussion above 
can be roughly summarized by the following 
\vskip0.1cm
\noindent
{\bf Algorithm inRS: NURBS-shaped indicator function}
\\ \\
{\bf input}: a point cloud $S$, the NURBS boundary parametrization $\Gamma(t)
=(\alpha(t),\beta(t))$  

\begin{itemize}

\item[$(i)$] determine the monotone boxes 
by computing the zeros of $\alpha'(t)$ and of $\beta'(t)$ as 
in (\ref{zerosalpha'beta'})-(\ref{mbox})  

\item[$(ii)$] compute a bounding box $B$ for $\Omega$ by the
extremal monotone boxes and collect the critical boundary points: 
vertical ($\alpha'=0$) or horizontal ($\beta'=0$) tangency points, 
$x$-turning or $y$-turning vertices  
 
\item[$(iii)$] for all $P\in S\cap B^c$: 
set $inRS(P)=0$
\item[$(iv)$] for all $P\in S\cap B$:
 
\begin{itemize} 
\item[$(v)$] compute $cross(P)$, the number of boundary 
intersections 
of a vertical ray from $P$ in the monotone boxes $\mathcal{B}_\ell$, 
with $\ell$ given in (\ref{whichboxes}) 
or (\ref{whichboxes2}) in case of box refinement (solving 
explicitly (\ref{intersections}) only if $P\in \mathcal{B}_\ell$) 

\item[$(vi)$] if the intersections do not include critical points 
(cf. Fig. \ref{fig_1CRI})\\ 
then $inRS(P)=cross(P)\;\mbox{mod}\,2$\\ else 
\begin{itemize}
\item[$(vii)$]
repeat step $(v)$ with an horizontal ray 

\item[$(viii)$] if the intersections do not include critical points\\ 
then $inRS(P)=cross(P)\;\mbox{mod}\,2$\\ else compute $inRS(P)=wind(P)$  
\end{itemize}
\end{itemize}
\end{itemize}
\noindent
{\bf output}: for all $P\in S$, $inRS(P)=1$ if $P\in S\cap \Omega$, 
$inRS(P)=0$ otherwise  
\vskip0.1cm

It is worth stressing some relevant features. 
First, we have constructed a suitable 
user-friendly structure for the NURBS boundary
parametrization, since it is apparently missing or 
at least complicated within basic Matlab. 

Moreover, if $P$ is a pointset 
instead of a single point, steps $(i)-(ii)$ can be clearly done once. 
Steps $(i)$, 
$(v)$ and $(vii)$ require solving polynomial equations, that can be 
conveniently done by the Matlab {\sc roots} function. Computation 
of $wind(P)$ in $(viii)$ has been implemented by Gaussian quadrature 
along the boundary but is more costly 
than $cross(P)$ (whenever the latter is feasible), so it has been 
reserved to dubious cases.    

In addition, we have tried optimizing 
all the algorithm blocks, conveniently managing the critical situations   
and implementing more features than those present in 
the polynomial spline version \cite{SV21-first} and in the rational splines 
alpha-version used in \cite{SV21}, that we call $inRS1$ below. 
For example, before resorting to 
$wind(P)$ we have provided the horizontal ray step after the vertical one, 
and we have added the possibility of knowing whether the point is 
in the domain interior or on the boundary (up to a suitable numerical 
tolerance).

The main improvement with respect to \cite{SV21-first} is however 
in the implementation of $(iv)$. Indeed, in \cite{SV21-first} the algorithm 
loops over the points in $S\cap B$ and for each point loops over the 
monotone boxes. With $M$ points and $N$ boxes the computational complexity    
is then $\mathcal{O}(MN)$.
 
Differently, in the present version of the algorithm we chose the 
strategy already used in \cite{E21} for linear polygons, based on 
point ordering and binary search, as follows: 

\begin{itemize}

\item orders the points in $S\cap B$ w.r.t. the $x$-variable  

\item loops over the monotone boxes

\begin{itemize}
\item for each monotone box finds all the points 
$P$ satisfying (\ref{whichboxes}) 
(or (\ref{whichboxes2})) by binary search followed by local linear search   

\item for all such points performs $(v)$ and possibly 
updates $cross(P)$ whenever a noncritical intersection is found 
\end{itemize}
eventually, for every point $P\in S\cap B$, 
either the overall number of crossings 
is computed and $(vi)$ can be done, or $P$ is recognized as ``dubious''

\item performs $(vii)$-$(viii)$ for all dubious points   

\end{itemize} 

Algorithm $inRS$ with the implementation just described will be 
called $inRS2$ below. 
With $M$ points and $N$ boxes the computational complexity is now 
reduced to  
$\mathcal{O}(M\log_2(M))$ for the ordering, plus $\mathcal{O}(N\log_2(M))$  
for the binary searches, plus $\mathcal{O}(M\nu)$ where $\nu$ is the mean 
number of monotone boxes to which a point in $S\cap B$ belongs (when 
the boxes are small and the points quasi-uniformly distributed 
this number, which corresponds to the number of equations 
like (\ref{intersections}) to be solved by {\sc roots}, 
can be very small).

\begin{figure}[!htbp]
  \centering
  {\includegraphics[scale=0.27,clip]{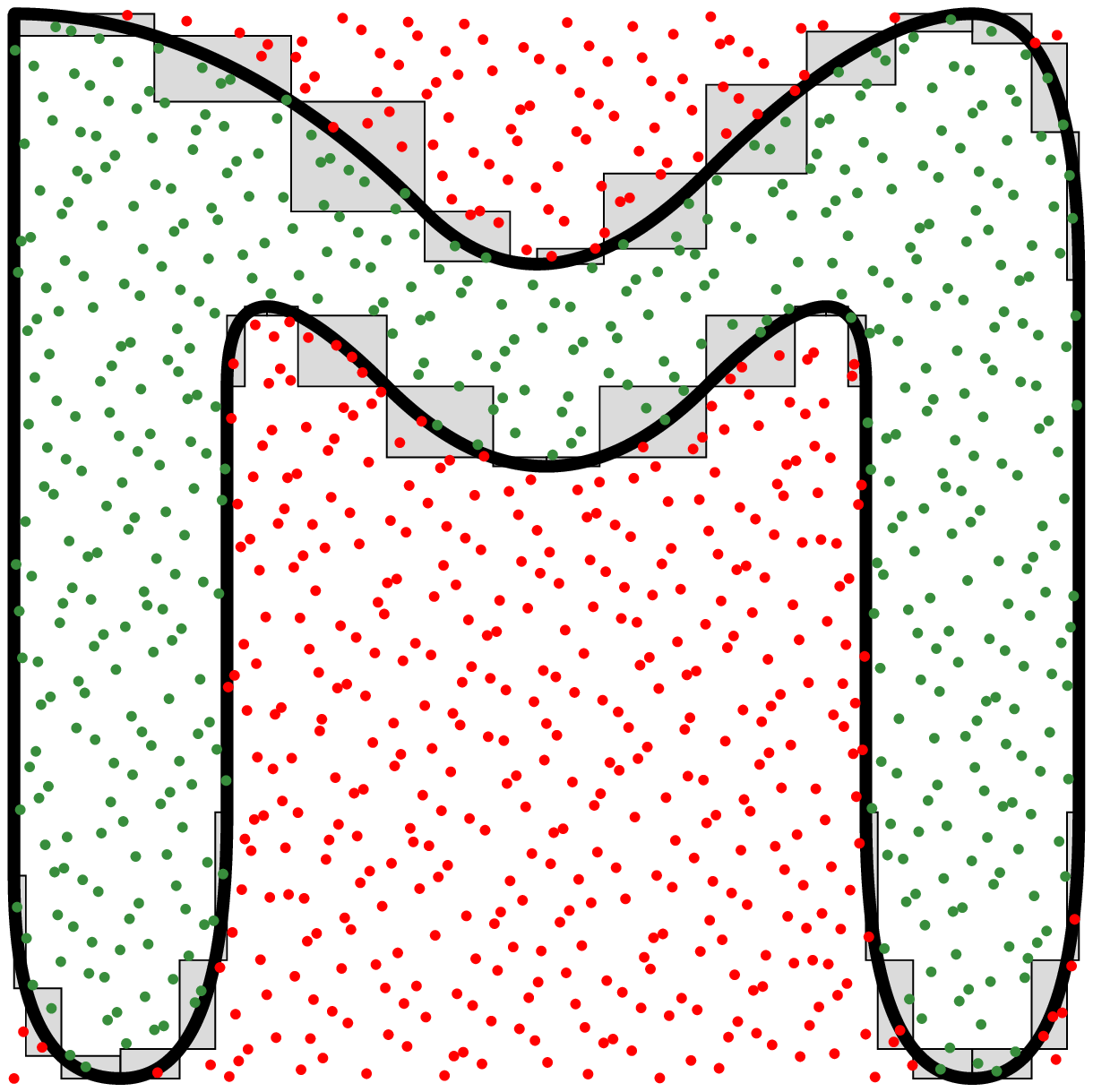}}
  {\includegraphics[scale=0.27,clip]{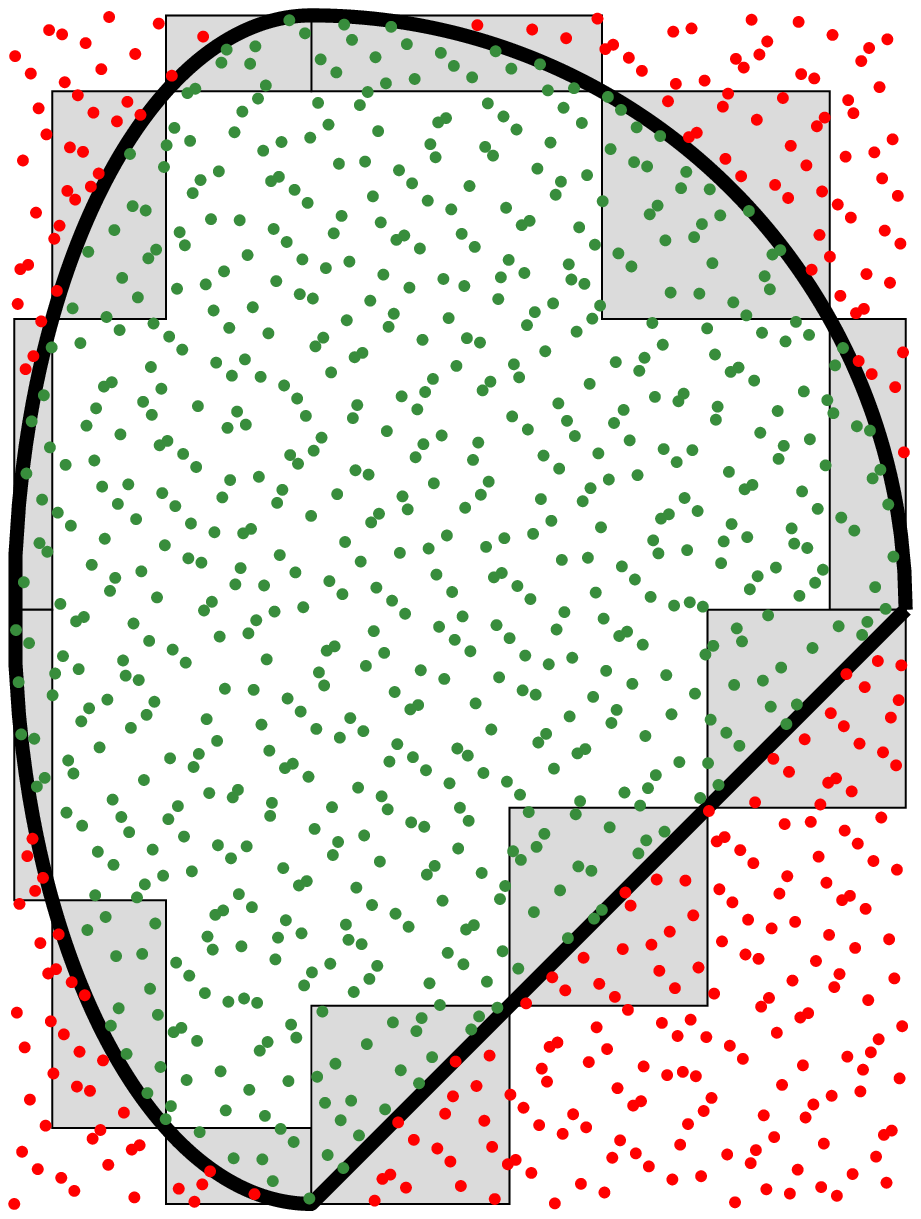}}
       {\includegraphics[scale=0.2,clip]{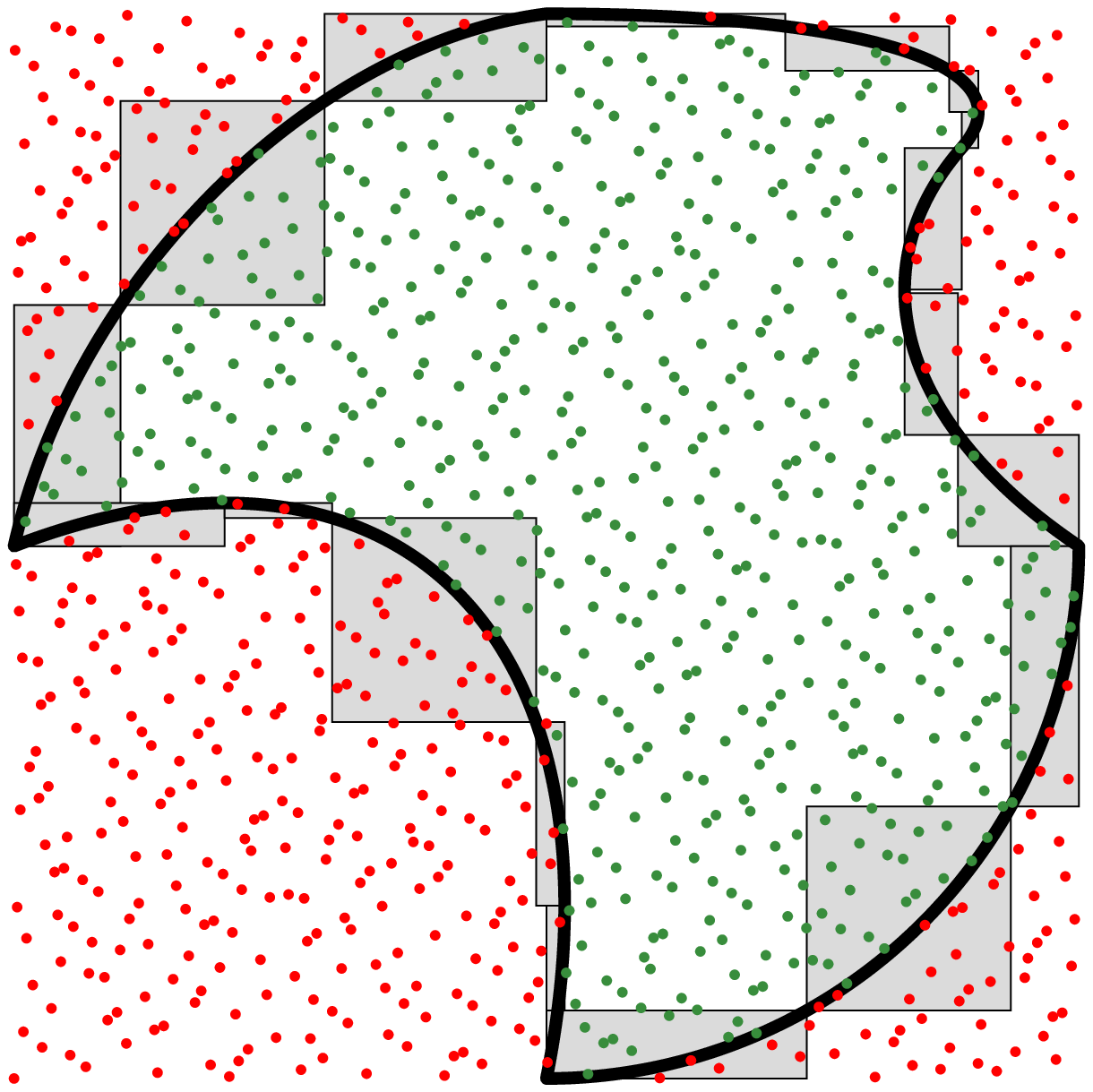}}
 \caption{Halton points inclusion for the three NURBS-shaped domains
of Fig. \ref{fig2}.}
 \label{fig4}
 \end{figure}

A graphical example concerning the inclusion check on 1000 Halton 
points of the domain bounding box is given in Fig. \ref{fig4}. 
We stress that equations like (\ref{intersections}) have to be solved only 
for the points belonging to the monotone boxes.     
In order to give an illustration of the algorithm performance, 
we display Table \ref{tab01} below, where we locate $10^i$ 
Halton points  
of the global bounding box, $i=3,4,5$, on the three Nurbs-shaped 
domains of Fig. \ref{fig2}; the speed-up is rounded 
to two significant figures and  
to manage CPU time fluctuations we have taken the median over 
100 runs of the algorithms.
As expected, the speed-up of $inRS2$ over $inRS1$ increases 
by increasing the cardinality, taking however into account 
that both have a fixed cost due to the construction of monotone boxes 
which is relevant at the lower cardinalities.  

Similar results are obtained with 
the same number of grid points or random points. The numerical tests  
of the present paper have been done on a M1-chip PC with 16 GB of RAM, 
running Matlab R2021b. All the Matlab codes and demos of the present version
are freely available at \cite{SV22-soft}.  

\begin{table}[!htbp]
\begin{center}\footnotesize
\begin{tabular}{|c | c| c | c | c |}
\hline
$\#$ & algorithm   & Fig. \ref{fig2}-left   & Fig. \ref{fig2}-center
& Fig. \ref{fig2}-right  \\
\hline
& $inRS1$  & 4.3$e-$03$s$  & 1.9$e-$03$s$ &  2.1$e-$03$s$\\
$10^3$ & $inRS2$ & 2.8$e-$03$s$  & 1.3$e-$03$s$  & 1.3$e-$03$s$\\
& speed-up & 1.5  & 1.5 & 1.6\\
\hline
\hline
& $inRS1$  & 1.8$e-$02$s$  & 8.8$e-$03$s$ & 9.5$e-$03$s$\\
$10^4$ & $inRS2$ & 4.4$e-$03$s$  & 3.2$e-$03$s$ & 3.2$e-$03$s$\\
& speed-up & 4.1  & 2.8 & 3.0\\
\hline
\hline
& $inRS1$  & 1.6$e-$01$s$  & 7.3$e-$02$s$ &  8.0$e-$02$s$\\
$10^ 5$ & $inRS2$ & 1.7$e-$02$s$  & 2.1$e-$02$s$ & 2.0$e-$02$s$\\
& speed-up & 9.4 & 3.5 & 4.0\\
\hline
\end{tabular}
\caption{\small{CPU time of
$inRS$ on the three 
NURBS-shaped domains of Fig. \ref{fig2}-\ref{fig3} with $\#$ Halton points 
of the corresponding bounding box; $inRS1$ is the alpha-version in 
\cite{SV21} whereas $inRS2$ is the present version.}}
\label{tab01}
\end{center}
\end{table}

\begin{Remark}
{\em 
It is worth recalling that the fact that all inequalities 
in the algorithm above, as usual in computational geometry, 
are checked up to 
a given tolerance, so that ultimately we only know whether a point $P$ 
is in or out a suitable neighborhood of the NURBS-shaped domain.    

Since there are clever implementations of the in-domain 
check for linear polygons, one may think that it could be better to 
approximate the NURBS-shaped boundary by a piecewise linear curve 
with a very high number of sides up to the given tolerance, and then 
to use such fast point-in-polygon algorithms (roughly, if a tolerance 
$\varepsilon$ is given, the number of sides 
is $O\left(\varepsilon^{-1/2}\right)$ 
since the error of linear approximation is $O(\Delta t ^2)$ 
for a piecewise $C^2$ NURBS parametrization).

In practice, however, for such small tolerances our implementation 
of the NURBS-shaped indicator function is faster, as it can be 
realized by Table \ref{tab02} below, where we compare the {\sc inpoly2} 
Matlab function \cite{E21,KKE21} (that largely overcomes 
the standard {\sc inpolygon} 
on a large number of trial points) applied to a polygonal approximation 
up to a $10^{-10}$ tolerance (order of $10^5$ sides), with our $inRS2$ 
on the three NURBS-shaped 
domains of Fig. \ref{fig2}-\ref{fig3}. 
Notice that the speed-up of 
$inRS2$ over {\sc inpoly2} tends to decrease by increasing the cardinality 
of the point set to be located. The dynamics of 
CPU times in the present 
range of cardinalities can be 
explained observing that, while that of $inRS2$ is substantially ruled by 
the number of points (the number of boxes being relatively small 
with respect to it) up to the fixed initial cost of the boxes construction, 
that of {\sc inpoly2} is instead substantially ruled 
by the huge number of polygon sides.    
\/}
\end{Remark}

\begin{table}[!htbp]
\begin{center}\footnotesize
\begin{tabular}{|c | c| c | c | c |}
\hline
$\#$ & algorithm   & Fig. \ref{fig2}-left   & Fig. \ref{fig2}-center
& Fig. \ref{fig2}-right  \\
\hline
& {\sc inpoly2}  & 2.1$e-$01$s$  & 5.2$e-$02$s$ &  6.5$e-$02$s$\\
$10^3$ & $inRS2$ & 2.8$e-$03$s$  & 1.3$e-$03$s$  & 1.3$e-$03$s$\\ 
& speed-up & 75  & 40 & 50\\
\hline
\hline
& {\sc inpoly2}  & 2.5$e-$01$s$  & 6.4$e-$02$s$ & 8.0$e-$02$s$\\
$10^4$ & $inRS2$ & 4.4$e-$03$s$  & 3.2$e-$03$s$ & 3.2$e-$03$s$\\  
& speed-up & 58 & 20 & 25\\
\hline
\hline
& {\sc inpoly2}  & 3.1$e-$01$s$  & 8.7$e-$02$s$ &  1.1$e-$01$s$\\
$10^ 5$ & $inRS2$ & 1.7$e-$02$s$  & 2.1$e-$02$s$ & 2.0$e-$02$s$\\
& speed-up & 18 & 4.1 & 5.5\\
\hline
\end{tabular}
\caption{\small{CPU time of  
the $inRS2$ algorithm on the three
NURBS-shaped domains of Fig. \ref{fig2}-\ref{fig3} with $\#$ Halton points
of the corresponding bounding box, compared with the Matlab 
{\sc inpoly2} function applied to an approximating 
polygon up to a $10^{-10}$ tolerance (order of $10^5$ sides).}}
\label{tab02}
\end{center}
\end{table}

\section*{Acknowledgements}

This work was partially supported by the DOR funds and the biennial
project BIRD 192932 of the University of Padova, and by the INdAM-GNCS, 
and has been accomplished within the RITA
Research ITalian network on Approximation and the UMI Group TAA
Approximation Theory and Applications.


\begin{thebibliography}{99}

\bibitem{BdVRV19}  L. Beir\~{a}o da Veiga, A. Russo and G. Vacca,
The Virtual Element Method with curved edges,
ESAIM Math. Model. Numer. Anal.,
53, 2019.

\bibitem{BDME16} 
L. Bittante, S. De Marchi and G. Elefante, 
A new quasi-Monte Carlo technique based on nonnegative least-squares 
and approximate Fekete points, Numer. Math. TMA., 9, 2016. 

\bibitem{DEBOOR}
C. de Boor,
\newblock A Practical Guide to Splines, Rev.ed.
Springer-Verlag, New York, 2001.

\bibitem{BPV20} L. Bos, F. Piazzon and M. Vianello,
Near G-optimal Tchakaloff designs,
Comput. Statistics, 35, 2020.








\bibitem{CB15} J.-S. Chen and T. Belytschko, Meshless and Meshfree Methods, 
in: B. Engquist Ed., Encyclopedia of Applied and Computational Mathematics, 
Springer, 2015, pp. 886--894. 


\bibitem{E21} D. Engwirda, INPOLY: A fast points-in-polygon test, 
GitHub, 2021: \url{https://github.com/dengwirda/inpoly}.  

\bibitem{FMC15} G.E. Fasshauer and M.J. McCourt, Kernel-based Approximation 
Methods using Matlab, Interdisciplinary Mathematical
Sciences, Vol. 19, World Scientific Publishing Co., Singapore, 2015.

\bibitem{GUND21}
D. Gunderman, K. Weiss and J.A. Evans,
Spectral mesh-free quadrature for planar regions bounded
by rational parametric curves,
Computer-Aided Design, 130, 2021.

\bibitem{H21}
S. Hayakawa, 
Monte Carlo cubature construction, 
Jpn. J. Ind. Appl. Math., 38, 2021.

\bibitem{H07} T.C. Hales, Jordan's Proof of the Jordan Curve Theorem, 
Studies in Logic, Grammar and Rethoric, 10 (23), 2007. 

\bibitem{HA01} K. Hormann and A. Agathos, The point in polygon problem 
for arbitrary polygons, Comput. Geom., 20, 2001. 

\bibitem{J1893} C. Jordan, Course d'analyse de l'\`{E}cole Polytechnique, 
Gauthier-Villars, Paris, 1893. 

\bibitem{KKE21} J. Kepner, A. Kipf, D. Engwirda, \& al., 
Fast Mapping onto Census Blocks, 
2020 IEEE High Performance Extreme Computing Conference (HPEC), 
arXiv:2005.03156v2. 

\bibitem{Kr99} S.G. Krantz, The Index or Winding Number of a Curve 
about a Point, \S 4.4.4 in Handbook of Complex Variables, 
Birkh\"{a}user, Boston,  1999.

\bibitem{PSV17} F. Piazzon, A. Sommariva and M. Vianello, 
Caratheodory-Tchakaloff Least Squares, Sampling Theory 
and Applications 2017, IEEE Xplore Digital Library, 
DOI: 10.1109/SAMPTA.2017.8024337. 

\bibitem{PIEGL}
L. Piegl,
\newblock The NURBS Book, Second Edition, Springer-Verlag,
Berlin-Heidelberg, 1997.

\bibitem{SFM11}
R. Sevilla and S. Fern\'andez-M\'endez,
Numerical integration over 2D NURBS-shaped domains with
applications to NURBS-enhanced FEM, 
Finite Elements in Analysis and Design,
47, 2011.


\bibitem{SV21-first}
A. Sommariva and M. Vianello,
Computing Tchakaloff-like cubature rules on spline curvilinear polygons, 
Dolomites Res. Notes Approx. DRNA, 14, 2021. 

\bibitem{SV21}
A. Sommariva and M. Vianello,
Low-cardinality Positive Interior cubature 
on NURBS-shaped domains,
preprint, 2021, available online at: 
\url{https://www.math.unipd.it/~marcov/pdf/nurbscatch.pdf}.

\bibitem{SV22-soft}
A. Sommariva and M. Vianello,
inRS: a Matlab code for the indicator function of NURBS-shaped planar 
domains, available online at:
\url{https://www.math.unipd.it/~alvise/software.html}.



\end{thebibliography}
\end{document}